  \documentclass[jaa]{degruyter-journal-a}      % Journal of Applied Analysis
\usepackage{mathrsfs}

\theoremstyle{plain} %default

\newtheorem{conj}{Conjecture}

\newtheorem{thm}{Theorem}

\theoremstyle{definition}

\theoremstyle{remark}

\numberwithin{equation}{section}

\DeclareMathOperator{\divgce}{Div}
\DeclareMathOperator{\pd}{\partial}

\DeclareFontFamily{OT1}{pzc}{}
\DeclareFontShape{OT1}{pzc}{m}{it}{<-> s * [1.10] pzcmi7t}{}
\DeclareMathAlphabet{\mathpzc}{OT1}{pzc}{m}{it}

\DeclareMathSymbol{\R}{\mathalpha}{AMSb}{"52}
\DeclareMathSymbol{\C}{\mathalpha}{AMSb}{"43}
\newcommand{\mbb}[1]{\mathbb{#1}}
\newcommand{\F}{\mbb{F}}

\newcommand{\bv}{\mathbf{v}}

\newcommand{\bw}{\mathbf{w}}

\newcommand{\fq}{\mathfrak{q}}

\title{Some variational principles  associated with ODEs of maximal symmetry. Part 1: Equations in canonical form}
\headlinetitle{Variational principles  associated with ODEs in canonical form}

\lastnameone{Ndogmo}
\firstnameone{Jean-Claude}
\nameshortone{J.\,C.~Ndogmo}

\addressone{Department of Mathematics and Applied Mathematics,
University of Venda,
P/B X5050, Thohoyandou,
Limpopo  0950}
\countryone{South Africa}
\emailone{ndogmoj@yahoo.com}

\abstract{Variational and divergence symmetries are studied  in this paper for linear equations of maximal symmetry in canonical form, and the associated first integrals are given in explicit form.  All the main results obtained  are formulated as theorems or conjectures for equations of a general order. Some of these results apply to linear equations of a general form and of arbitrary orders or having a symmetry algebra of arbitrary dimension.}



\researchsupported{This research is supported by  research grants from the University of Venda and from the NRF of South Africa.}
\begin{document}
%%%%%%%%%%%%%%%%%%%%%%%%%%%%%%%%

%%%%%%%%%%%%%%%%%%%%%%%%%%%%%%%%

\section{Introduction}
\label{s:intro}
%%%%%%%%%%%%%%%%%%%%%%%%%%
%%%%%%%%%%%%%%%%%%%%%%%%

Due to some  recent results on linear ordinary differential equations ({\sc lode}s) of  maximal symmetry \cite{KM, ML, JF, J14}, it should be possible by now to describe in their most general form the variational principles associated with this class of equations and with the whole class of nonlinear  {\small \sc ode}s of maximal symmetry, that is, having a Lie point symmetry algebra of maximal dimension. In particular, it should be possible to obtain  the precise conditions for the existence of variational symmetries, and the related first integrals thereof. \par

Pursuing a work undertaken by Lie \cite{LieCanonic}, Krause and Michel \cite{KM} have indeed given a specific proof of the fact that a {\small \sc lode} has maximal symmetry if and only if it can be reduced by a point transformation to the form $y^{(n)}=0,$ which will henceforth be referred to as the canonical form. In the same theorem of \cite{KM} proving a special case of a theorem of \cite{LieCanonic}, the {\small \sc lode} is also shown to have maximal symmetry if and only if it is iterative. Mahomed and Leach \cite{ML} also gave a partial characterization of this class of equations, and an expression for the coefficients of the linear equation in terms of the coefficients of the term of third highest order,  for orders $k$ of the equation such that $3 \leq k \leq 8.$ A general expression for the coefficients of these {\small \sc lode}s was recently obtained in terms of the parameters of the source equation, as well as the operator generating an equation of maximal symmetry of a general order, together with several important transformation and other properties of these equations \cite{JF, J14}.\par

It is well known that a {\small \sc lode} can always be reduced by a point transformation to the normal form
\begin{equation} \label{nor1}
y^{(n)} + A_n^2 y^{(n-2)}+ \dots + A_n^{n-1} y^{(1)} + A_n^{n} y =0,
\end{equation}
in which the coefficient of the term of second highest order has vanished. A specific property satisfied by {\small \sc lode}s of maximal symmetry in the form \eqref{nor1} is that $n$ linearly independent solutions $s_k$ can be obtained in the form
\begin{equation}\label{sk}
s_k=  u^{n-k-1} v^k, \qquad 0 \leq k \leq n-1, \qquad n \geq 2,
\end{equation}
for some functions $u$ and $v$ of $x.$ Moreover, when such equations are in normal form \eqref{nor1}, the Wronskian  $\mathpzc{w}(u,v)= \det {\small \begin{pmatrix} u & v\\ u'  & v' \end{pmatrix} }$  is a nonzero constant and will be normalized to one. It also follows from \eqref{sk} that $u$ and $v$ are two linearly independent solutions of the second-order source equation
\begin{equation} \label{srce}
y_{xx} + \fq(x)y=0,\qquad \fq= A_2^2.
\end{equation}
Using in particular this fact, the symmetry generators of equations of the form \eqref{nor1} of maximal symmetry were all obtained in \cite{KM} in terms of two linearly independent solutions $u$ and $v$  of \eqref{srce}. In the sequel, when we refer to \eqref{nor1}, we shall assume that it is  the most general form of a {\small \sc lode} of maximal symmetry, unless otherwise stated.\par

 In this paper, by exploiting the expression of the symmetry generators for equations of the form \eqref{nor1} of maximal symmetry, we investigate for the corresponding class of equations in canonical form $y^{(n)}=0$, the existence, and carry out the determination of  variational and  divergence symmetry algebras, and first integrals in the sense given in \cite{olv1}. The results thus obtained ultimately apply to {\small \sc lode} of a general order, and often to any {\small \sc lode} of the most general form regardless of the dimension of its symmetry algebra.  These results  are  extended in Part 2 \cite{part2} of this series of papers to the most general form \eqref{nor1}, and then to the whole class of nonlinear  {\small \sc ode}s of maximal symmetry. The field $\F$  of scalars will be assumed to be the field $\R$ of real numbers, although $\F$ could equally be taken to be $\C$.

\section{Some relevant facts about variational symmetries}
\label{s:basic}
Let $G$ be a local group of transformations acting on the space $X \times Y  \subset  \F^{\,p} \times \F^q,$ where $X$ represents the space of independent variables and $Y$ the space of dependent variables. Denote by $f[y]= f(x, y^{(n)})$ a differential function depending on $x,\, y= y(x),$ and the derivatives of $y$ up to an arbitrary but fixed order $n.$ \par

Given an open and connected subset $\Omega \subset X,$ with smooth boundary $\partial \Omega,$ the problem of finding the extremals of the functional
\begin{equation}\label{varprob}
\vartheta[y] = \int_{\Omega} L(x, y^{(n)}) dx
\end{equation}
is called a variational problem and the integrand $L= L(x, y^{(n)})$ is referred to as the Lagrangian. Smooth extremals of \eqref{varprob} are attained on functions where the variational derivative $\delta \vartheta  [y]$ of the functional $\vartheta  [y]$ vanishes. The variational derivative $\delta \vartheta  [y]$ has a very simple expression in terms of the Euler operator and $L.$ More succinctly, every smooth extremal $y= f(x)$ of \eqref{varprob} must be a solution of the associated Euler-Lagrange equations
\begin{subequations}\label{EL}
\begin{align}
E_\alpha (L)&=0, \qquad \alpha= 1, \dots, q \label{EL1}\\
\intertext{where }
E_\alpha &= \sum_J (-D)_J \frac{\partial}{\partial y_J^\alpha} \label{EL2}
\end{align}
\end{subequations}
is the $\alpha$-th Euler operator. In \eqref{EL2}, $D$ is the total differential operator and $J= (j_1, \dots, j_k),$ with $1 \leq j_i\leq p.$

In the actual case of ordinary differential equations that we shall be dealing with, we have $p=q=1.$ By letting $D_x$ denote the total differential operator with respect to the single variable $x,$  the Euler-Lagrange equations take the form
\[
0= E(L) =  \sum_{j=0}^\infty (-D_x)^j \frac{\partial }{\partial y_j}=  \frac{\partial L}{\partial y} - D_x \frac{\partial L}{\partial y_x}  + \dots + (-1)^n D_x^n \frac{\partial L}{\partial y_n},
\]

where $y_j= d^j y/ d x^j,$ and $y_x, y_{x x}, \dots$ represent, respectively, the derivatives  $d y/ d x, d^2 y/ d x^2, \dots.$
This is a $2n$-th order {\sc ode} provided $L$ satisfies the nondegeneracy condition $\partial ^2 L/ \partial y_n^2 \neq 0.$\par

The local group of transformations $G$ acting on an open subset $M \subset \Omega \times Y$ is called a variational symmetry group of \eqref{varprob} if it leaves the corresponding integral invariant. The precise definition of  a variational symmetry group is a bit too detailed for our discussion and can be found in
\cite[P 157]{olv1} (see also \cite{blum}). The infinitesimal version of this definition is however fairly simple to express. Indeed, let
\begin{equation}\label{v}
\mathbf{v} = \sum_{i=1}^p \xi^i(x,y) \frac{\partial }{\partial x^i} + \sum_{\alpha=1}^q \psi_\alpha (x,y) \frac{\partial }{\partial y^\alpha}
\end{equation}
be an infinitesimal generator of $G,$ and for every  $p$-tuple $B=(B_1, \dots, B_p)$ of smooth functions on $\Omega,$ denote by $\divgce B = \sum_{i=1}^p D_{x^i} B_i$ the total divergence of $B.$ We shall often use the notation $y_{(m)}$ to denote $y$ and all its derivatives up to the order $m,$ and the $m$th prolongation of a vector field $\bv$ will be denoted by $\bv^{[m]}.$   With this notation,  the connected transformation group $G$ acting on $M \subset \Omega \times Y$ is a variational symmetry group of \eqref{varprob} if and only if
\begin{equation} \label{vsym}
 \bv^{[n]} (L) + L \divgce \xi =0
\end{equation}
%%%%%%
for all $(x, y_{(n)})$ in the extended jet space $M^{(n)}$ of $M,$ and for every infinitesimal generator $\mathbf{v}$ of the form \eqref{v} of $G.$ A \emph{conservation law} for the Euler-Lagrange equations \eqref{EL1} is an expression $\divgce F=0$ for some $p$-tuple $F= (F_1,\dots, F_p),$  which vanishes identically on solutions of \eqref{EL}. Such an expression will usually be identified with the divergence expression $\divgce F$ on the left hand side of the equation.\par

Let
\begin{subequations}\label{onvq}
\begin{align}
\mathbf{v}^Q &= \sum_{\alpha=1}^q Q_{\alpha} \frac{\partial }{\partial y^\alpha} \label{vq}\\
\intertext{be the characteristic form of $\mathbf{v}$ as given in \eqref{v}. Thus}
Q_\alpha &= Q_\alpha (x,y^{(1)}) = \psi_\alpha - \sum_{i=1}^p \xi^i y_i^\alpha,\quad y_i^\alpha = \partial y^\alpha / \partial x^i, \label{vqfrmla}
\end{align}
\end{subequations}

 and the vector $Q= (Q_1, \dots, Q_q)$ is the corresponding  characteristic for $\mathbf{v}.$ Noether's Theorem asserts that for every infinitesimal variational symmetry $\mathbf{v}$ as given in \eqref{v} its characteristic $Q$ is also the characteristic of a conservation law for the Euler-Lagrange equations $E(L)=0.$ That is,
 $$
 \divgce F= Q\cdot E(L) = \sum_{\alpha=1}^q Q_\alpha E_\alpha (L)
 $$
 is a conservation law in characteristic form, and for some $p$-tuple $F=(F_1, \dots, F_p)$ of smooth functions on $\Omega.$\par

 On the other hand, a \emph{divergence symmetry} of the functional $\vartheta[y]= \int_\Omega L dx$ is a vector field of the form \eqref{v} which satisfies
 \begin{equation} \label{dsym}
 \mathbf{v}^{[n]} (L) + L \divgce \xi = \divgce B
 \end{equation}
for all $x,y \in M$ and for some $p$-tuple $B=B(x, y^{(m)})= (B_1, \dots, B_p)$ of smooth functions on $M$. Thus the divergence symmetries of \eqref{varprob} satisfy a more relaxed condition of the form \eqref{vsym} required for variational symmetries. It turns out that both variational symmetries and divergence symmetries of the variational problem \eqref{varprob} are all symmetries of the associated Euler-Lagrange (EL) equations. Moreover, the characteristic of each divergence symmetry is also the characteristic of a conservation law of the associated EL equations.\par

Denote by $\mathcal{A}$ the space of smooth differential functions $F(x, y^{(n)})$ on $M^{(n)} \subset X \times Y^{(n)},$ and by $\mathcal{A}^m$ the vector space of $m$-tuples of differential functions $F[y]= (F_1[y], \dots, F_m[y]),$ where each $F_j \in \mathcal{A}.$  It is well known that an ordinary vector field of the form \eqref{v} is a divergence symmetry if and only if its evolutionary form $\mathbf{v}^Q$ is.\par

For $F\in \mathcal{A}^r$ denote by $D_F$ its Fr\'echet derivative and by $D_F^*$ its adjoint Fr\'echet derivative, i.e. the adjoint operator of $D_F.$ Then it follows from Leibnitz rule that
\begin{equation} \label{mplier}
E (F \cdot Q) = D_F^* (Q) + D_Q^* (F), \quad \text{ for all } F, Q \in \mathcal{A}^r.
\end{equation}
Moreover, a  $q$-tuple $Q \in \mathcal{A}^q$ is the characteristic of a conservation law associated with a given system $(\Delta=0)$ of differential equations if and only if
\begin{equation} \label{adjcond}
D_\Delta^* (Q) +D_Q^* (\Delta) =0
\end{equation}
for all $(x,y) \in M.$ We shall usually  make use of this condition to identify divergence symmetries of the equations we wish to investigate. This means that divergence symmetries may be defined even for equations which do not admit a Lagrangian formulation. In this way, the definition of divergence symmetries naturally extends in particular to equations of any (even or odd) order.\par

In fact \eqref{mplier} and \eqref{adjcond} are the justification of so-called multiplier/ integrating factor method \cite{anco, magan-mahomed, bashingwa-kara, delarosa} for finding conservation laws  associated with a system of differential equations regardless of its symmetry properties and its equivalence or not to a Lagrangian equation.

%%%%%%%

\section{Lagrangians}\label{s:lagrang}

Since variational symmetries and divergence symmetries of the variational problem $\vartheta[y] = \int L dx$ are also symmetries of the associated EL equations, one way to find them is to identify in the Lie point symmetry algebra of the equations the most general linear combination of the symmetry generators that satisfy the required conditions. However, it follows from condition \eqref{vsym} that for the identification of variational symmetries, we shall need to know the associated Lagrangian function. \par

First of all, it can readily be verified that all {\small \sc lode}s of maximal symmetry (and necessarily of even orders) up to the order $8$ satisfy Helmholtz conditions for the inverse problem of the calculus of variations. That is, these equations, obtained in \cite{ML} and \cite{J14}, are the EL equations  of some variational problem. Their explicit expressions  are obtained  as follows for each order $n=2,4,6,8$ in terms of the single arbitrary coefficient $\fq= \fq(x).$
\begin{subequations}\label{o}
\begin{align}
0=& \fq y + y_{xx}\label{o2}\\
0=&10 y_x \fq_x+10 \fq y_{xx}+3 y \left(3 \fq^2+\fq_{x x}\right)+y^{(4)} \label{o4}\\[1.7mm]
 \begin{split}0=&7 y_{xx} \left(37 \fq^2+9 \fq_{x x}\right)+70 \fq_x y^{(3)}+y_x \left(518 \fq \fq_x+28 \fq^{(3)}\right)\\
 &+35 \fq y^{(4)}+5 y \left(45 \fq^3+26 \fq_x^2+31 \fq \fq_{x x}+\fq^{(4)}\right)+y^{(6)}
 \end{split}\label{o6}\\[1.7mm]
\begin{split} 0=&168 y^{(3)} \left(47 \fq \fq_x+ 2 \fq^{(3)}\right)+42 \left(47 \fq^2+9 \fq_{x x}\right) y^{(4)}\\
& +4 y_{xx} \left(3229 \fq^3+1773 \fq \fq_{x x}+45 \left(33 \fq_x^2+\fq^{(4)}\right)\right)\\
& +6 y_x \left(6458 \fq^2 \fq_x+524 \fq \fq^{(3)}+9 \left(132 \fq_x \fq_{x x}+\fq^{(5)}\right)\right) \\
& +252 \fq_x y^{(5)}+84\fq y^{(6)} +7 y \bigg(1575 \fq^4+1654 \fq^2 \fq_{x x}+153 \fq_{2}^2\\
&+226 \fq_x \fq^{(3)}+8\fq \left(347 \fq_x^2+10 \fq^{(4)}\right)+\fq^{(6)}\bigg) +y^{(8)}\label{o8}
\end{split}
\end{align}
\end{subequations}

Although calculations will be explicitly performed for all equations of maximal symmetry of order $n$ with $3\leq n \leq 8,$ for the sake of conciseness we do not list odd order equations. However, using a  recent result of \cite[Theorem 3.1]{JF}, these equations can be easily  generated for all orders. In the sequel, unless otherwise specified the notation $\Delta_n[y]=0$ represents a general  $n$th order {\small \sc lode} of maximal symmetry.\par

For any linear $n$th order equations $\Delta_n[y]=0,$  we can always take the associated Lagrangian to be $\frac{1}{2}y \Delta_n[y].$ However, the latter expression is of the same (even) order $n$ as the equation, and thus it can be simplified to yield an expression of order $n/2,$ by subtracting from it null Lagrangians of order higher than $n/2.$ To perform such a reduction, it suffices to note that each term in the expression $\frac{1}{2} y \Delta_n [y]$ has the form $f(x) y^{(s)}$ for a certain function $f(x).$ Each such term may be put into the form

\begin{equation} \label{fys}
f(x) y^{(s)} = D\left[  \sum_{j=0}^k (-1)^k D^k (f) y^{s-k-1} \right] + (-1)^k D^k (f) y^{s-k},
\end{equation}
for some $0 \leq k <s.$ This clearly reduces the corresponding term in the expression of the Lagrangian to $(-1)^k D^k y ^{(s-k)}.$ Thus each term $f(x) y^{(s)}$ of order $ s > n/2$ can be reduced to a term of order $n/2$ by letting $k= s- n/2.$  \par

Proceeding in this way, Lagrangians $L_n$ of order $n/2$ for  {\small \sc lode}s \eqref{nor1} of maximal symmetry  are readily obtained and those of order $n=2, 4, 6, 8$ are given as follows.
\begin{subequations} \label{L}
\begin{align}
L_2 &=  \frac{1}{2} \left[ - y_x^2 + \fq y^2\right]  \label{L2}  \\
L_4 &= \frac{1}{2} \left[ y_{x x}^2 - 10 \fq y_x^2 + 3 (3 \fq^2 + \fq_{x x} ) y^2    \right]  \label{L4}  \\
\begin{split}L_6 &= \frac{1}{2} \bigg[7 y y_{x x} \left(37 \fq^2+9 \fq_{x x}\right)+5 y^2 \left(45 \fq^3+26 \fq_x^2+31 \fq \fq_{x x}+\fq^{(4)}\right)\\
&\quad +y y_x \left(518 \fq \fq_x+28 \fq^{(3)}\right) -35 \fq y_x y^{(3)}+35 y \fq_x y^{(3)}-y_{3}^2\bigg]
\end{split} \label{L6}\\
\begin{split} L_8 &=  \frac{1}{2} \bigg[168 y y^{(3)} \left(47 \fq \fq_x+2 \fq^{(3)}\right)-84 y_x \fq_x y^{(4)}+84 \fq
y_{x x} y^{(4)}\\
&\quad  + 4 y y_{x x} \left(3229 \fq^3+1773 \fq \fq_{x x}+45
\left(33 \fq_x^2+\fq^{(4)}\right)\right) \\
&\quad +6 y y_x \left(6458 \fq^2 \fq_x+524 \fq \fq^{(3)}+9
\left(132 \fq_x \fq_{x x}+\fq^{(5)}\right)\right)\\
&\quad +42 y \left(47 \fq^2+5 \fq_{x x}\right) y^{(4)} +7 y^2 \bigg(1575 \fq^4+1654 \fq^2 \fq_{x x}+153 \fq_{2}^2\\
&\quad +226 \fq_x \fq^{(3)}+8 \fq \left(347 \fq_x^2+10 \fq^{(4)}\right)+\fq^{(6)}\bigg) +y_{4}^2 \bigg] \label{L8}
\end{split}
\end{align}
\end{subequations}

\section{Equations in canonical form }\label{s:canonic}

In order to gain a better insight into the properties of variational symmetries of {\small \sc lode}s  of maximal symmetry, it would be more convenient to look at the much simpler case of equations in canonical form $y^{(n)}=0.$ This is also useful for the extension of the results to the more general form \eqref{nor1}, and to the entire class of nonlinear {\small \sc ode}s of maximal symmetry as we shall do in Part 2 \cite{part2} of this series of two papers. \par

It was systematically shown in \cite{JF} that in its normal form \eqref{nor1} the most general {\small \sc lode} of maximal symmetry depends on a single arbitrary function which can be taken to be the coefficient $\fq=\fq(x)$ of the term of third highest order.  Moreover, the general expressions for these equations listed in \eqref{o} shows that choosing the canonical form $y^{(n)}=0$ amounts to setting $\fq=0$ in the general expression of the equation.\par

The expression \eqref{o2} for {\small \sc lode}s of second-order is not so special in a certain sense, because all {\small \sc lode}s of order 2 can be reduced by a point transformation to the free fall equation $y_{xx}=0.$  This is not the case however for equations of order $n\geq 3,$ and we shall thus restrict our attention to such equations. On the other hand, the symmetry algebra of such {\small \sc lode}s \eqref{nor1} of maximal symmetry and of the most general form is well known for an arbitrary order $n$ \cite{KM}.  First of all, by maximality the symmetry algebra has dimension $n+4$ for each equation of order $n.$ The infinitesimal generators are given by
\begin{subequations} \label{ig}
\begin{align}
V_k &= s_k \pd_y = u^{n-(k+1)} v^k \pd_y, \qquad 0\leq k \leq n-1 \label{ig1}\\
W_y &= y \pd_y    \label{ig2}\\
F_n &= u^2 \pd_x + (n-1) u u' y \pd_y    \label{ig3}\\
G_n &= 2 u v \pd_x + (n-1) (u v' + u' v) y \pd_y  \label{ig4}\\
H_n &= - v^2 \pd_x - (n-1) v v' y \pd_y,    \label{ig5}
\end{align}
\end{subequations}

where the $s_k$ are as in \eqref{sk} $n$ linearly independent solutions of the $n$th order equation, while $u$ and $v$ are two linearly independent solutions of the second order source equation \eqref{srce}.\par

It should also be mentioned here that the symmetries $V_k$  are often referred to as  solution symmetries and they generate the $n$-dimensional abelian Lie algebra $\mathcal{A}_n,$ while $W_y$ is often called the homogeneity symmetry. On the other hand, the Lie algebra $\mathfrak{s}$ generated by $F_n, G_n$ and $H_n$ is isomorphic to $\mathfrak{sl}_2 \equiv \mathfrak{sl}(2,\F)$ and the full symmetry algebra $\mathfrak{g}_n$ of \eqref{nor1} is the semi-direct sum $\mathfrak{g}_n= (\mathcal{A}_n \dotplus \F W_y) \ltimes \mathfrak{s},$ where the operator $\ltimes$ denotes a semi-direct sum of Lie algebras. We shall often denote by $\langle X_1,\dots,X_r \rangle$ a Lie algebra generated by the  vectors $X_1,\dots X_r.$ Thus for instance $\mathcal{A}_n= \langle V_0, \dots, V_{n-1} \rangle. $

The determination of \emph{first integrals} of linear equations of the form $y^{(n)}=0$ was carried out in \cite{flessas} with a different method not based on variational principles, but on a direct method,  by solving certain Lagrange equations to find the first integrals. Moreover, the focus on that paper was on the expression of first integrals in terms of functionally independent ones, for equations of order not exceeding $6.$\par

 In this section, we are interested in exhibiting the most general functional relation between the first integrals and the variables and parameters upon which they depend, for equations in canonical form of a general order. We also determine for such equations their divergence symmetry algebras, as well as their variational symmetry algebras in the case of equations admitting a Lagrangian description. Although calculations are performed most often only for equations of low orders not exceeding $8$, all the main results obtained are extended to {\small \sc lode}s of arbitrary orders.\par

For  a given Lagrangian $L$  corresponding to an {\sc ode} $\Delta_n[y]\equiv \Delta=0,$ set
\begin{align}
\mathscr{S}(\mathbf{v}) &=  \mathbf{v}^{[n]} (L) + L \divgce \xi,\qquad (\text{ for $n$ even}) \label{plv} \\
\mathscr{D}(\mathbf{v}) &= D_\Delta^* (Q) + D_Q^* (\Delta), \label{adjc}
\end{align}

where $\mathbf{v}= \xi(x,y) \pd_x + \psi(x,y)\pd_y$ has characteristic $Q.$ Then the operators $\mathscr{S}$ and $\mathscr{D}$ are clearly linear functions of $\mathbf{v},$ and $\mathscr{D}(\mathbf{v})$ is defined even for non-Lagrangian equations. It is also clear from \eqref{vsym} and \eqref{adjcond} that a (Lie point) symmetry $\bv$ of $\Delta=0$ is a variational symmetry if and only if
\begin{subequations}\label{v&d-cdt}
\begin{align}
\mathscr{S}(\mathbf{v})=& 0, \label{vsymeq} \\
\intertext{and a divergence symmetry   if and only if}
\mathscr{D}(\mathbf{v}) =&0  \label{dsymeq}.
\end{align}
\end{subequations}

For each order $n \geq 3,$ an arbitrary symmetry vector has the general form

\begin{subequations} \label{eq:bw}
\begin{align}
\mathbf{w} \equiv \mathbf{w}^n&=  \mathbf{w}_1 + \mathbf{w}_2  + \lambda W_y, \label{bw0} \\[-2mm]
\intertext{where \vspace{-3mm}}
\mathbf{w}_1 \equiv \mathbf{w}_1^n&= \sum_{k=0}^{n-1} a_k V_k,\quad \text{ and }\quad
\mathbf{w}_2 \equiv \mathbf{w}_2^n&= \alpha F_n+ \beta G_n + \gamma H_n,
\end{align}
\end{subequations}
for some arbitrary constants $a_k$ and $\alpha, \beta$ and $\gamma.$ In the actual case where $\fq=0,$ we shall let $u=1$ and $v=x$ be the two linearly independent solutions of the source equation \eqref{srce}. Let
\begin{equation} \label{Q_*^n}
F_n^Q= Q_F^n \pd_y,\quad G_n^Q= Q_G^n \pd_y,\quad \text{ and }H_n^Q= Q_H^n \pd_y.
\end{equation}
It follows from \eqref{onvq} that with the chosen values of $u=1$ and $v=x,$ the expressions for the characteristics are fairly simple, and are reduced to
\begin{align}
V_k^Q &= x^k \pd_y,\quad  \notag\\
Q_F^n &= - y_x,\quad Q_G^n = (n-1) y - 2 x y_x,\quad Q_H^n= (1-n) x y + x^2 y_x. \label{vqu=1}
\end{align}
Note that notation for expressions involving the variables $\fq, u,$ or $v$ will be preserved even when these variables are assigned the particular values $\fq=0, u=1$ or $v=x.$
%%%%%%%%%%%%%%%%%%%%%%%%%%%%%%%%%%%%%
\begin{thm}\label{th:wy}
For $n\geq 2,$ the homogeneity symmetry $W_y= y \pd_y$ is a divergence symmetry of any given $n$th-order {\small \sc lode} $\Delta_n[y]=0$  of maximal symmetry if and only if $n$ is odd. In particular for $n$ even, $W_y$ is not a variational symmetry.
\end{thm}
\begin{proof}
Since the characteristic $Q$ of   $W_y$ is $y,$  and $D_y^*$ is the identity operator, it follows from the fact that  the linear equations  $\Delta_n[y]=0$ are self-adjoint for $n$ even and skew-adjoint for $n$ odd that $\mathscr{D} (W_y)=0$ for $n$ odd and $\mathscr{D} (W_y)= 2 \Delta_n[y]$ for $n$ even. The first part of the theorem  then follows, and the second part is due to the fact that any variational symmetry is also a divergence symmetry. \par

The other way to see more directly in the case of Lagrangian equations why $W_y$ may not be a variational symmetry is to notice that in the expression of the Lagrangian $L_n$ obtained by reducing $\frac{1}{2} \Delta_n[y]$ to a  Lagrangian of order $n/2$ using  formula \eqref{fys}, the term $(-1)^{n/2} (D^{n/2} y)^2$ corresponding to the reduction of the term $\frac{1}{2} y y^{(n)}$ in $\frac{1}{2} y \Delta_n[y]$ is the only one with a square derivative of highest order. Moreover, $L_n$ is clearly a quadratic differential function of $y,$ and $W_y^{[n/2]}= \sum_{j=0}^{n/2} y_j \pd_{y_j}.$ Consequently, since in the variational symmetry expression \eqref{plv} we have $\xi=0,$ $\mathscr{S}(W_y)$ always contains the nonzero term $(-1)^{n/2} (D^{n/2} y)^2$ as the only one with a square derivative of highest order, and hence $\mathscr{S}(W_y)\neq0.$
\end{proof}

\begin{thm} \label{th:solsym2}
For all $n \geq 2,$ the Lie  algebra of divergence symmetries generated by the solutions symmetries $V_k$ of \eqref{nor1} is the whole of $\mathcal{A}_n= \langle V_0, \dots, V_{n-1} \rangle.$ In other words, $\mathcal{A}_n$ is always a subalgebra of the divergence symmetry algebra of $\Delta_n[y]=0.$
\end{thm}
\begin{proof}
For $\mathbf{v}= V_k= s_k \pd_y,$ since the characteristic $Q=Q(x)$ of $\mathbf{v}$ is a solution to the equation, we have $D_\Delta^*(Q) = \Delta (Q) =0.$ Moreover, $D_Q^*=0$ in this case where $Q$ depends only on $x$  and thus $\mathscr{D} (\mathbf{v})$ is always equal to zero. This completes the proof of the theorem.
\end{proof}

The proof of Theorem \ref{th:solsym2} shows that it is in fact valid for all linear equations of any order, and not just for those of maximal symmetry. It follows from Theorem \ref{th:wy} and Theorem \ref{th:solsym2} that in practice, to find the Lie algebra of divergence symmetries from the full Lie point symmetry algebra it will suffice to choose the vector $\bw$ in \eqref{bw0} in the form $\bw = \bw_2$ for odd order equations, and $\bw= \bw_2 + \lambda W_y$ for even order equations. Once the Lie subalgebra of divergence symmetries has been found, in order to find the corresponding subalgebra of variational symmetries the vector $\bw$ to be tested will then be chosen as the most general vector in the Lie subalgebra of divergence symmetries just found.

\subsection{Even order equations}
%%%%%%%%%%%%%%%%%%%%%%%%%%
 It follows from the above remark that for even order equations $y^{(n)}=0,$  the general vector $\bw$ of the form \eqref{bw0}  should be chosen as
 \[ \bw= \bw_2 + \lambda W_y = \alpha F_n + \beta G_n + \gamma H_n +  \lambda W_y.\]

 On the other hand, using the expressions of the characteristics  in \eqref{vqu=1},  a straightforward computation shows that
\[ \mathscr{D} (\mathbf{w}) \equiv \mathscr{D}(\alpha F_n + \beta G_n + \gamma H_n +  \lambda W_y)= 2 \lambda y^{(n)}.\]
 Therefore, due to the linearity of the operator $\mathscr{D}$ the most general divergence symmetry generator of $y^{(n)}=0$ is given by
\begin{equation} \label{dsym40}
\mathbf{w}= \alpha F_n + \beta G_n + \gamma H_n+\sum_{k=0}^{n-1} a_k V_k, \qquad \text{ for } n=4,6,8.
\end{equation}
In other words, among the Lie point symmetry group generators only $W_y= y \pd_y$ is not a term in the expression of the most general divergence symmetry vector.\par

Since for equations in canonical form $y^{(n)}=0$ we have $\fq=0$ in \eqref{o}, the expression of Lagrangians for equations of even orders $n=4,6,8$ reduces to the simpler form
\begin{equation} \label{Ln0}
L_n= \frac{(-1)^{n/2}}{2} \left(D^{n/2} y \right)^2.
\end{equation}
In fact, one can easily verify that this expression for $L_n$ is valid for all orders of equations in canonical form. The evaluation of $\mathscr{S}(\mathbf{w})$   is thus fairly simple in this case, and we let an arbitrary variational symmetry vector of $y^{(n)}=0$ be of the form \eqref{dsym40}. A direct computation of $\mathscr{S} (\mathbf{w})$ shows that
\begin{align*}
\text{For $n=4,$ }\quad \mathscr{S} (\mathbf{w}) &= 2(a_2 + 3 a_3 x - 2 \gamma y_x) y_{xx}\\
\text{For $n=6,$ }\quad \mathscr{S} (\mathbf{w}) &= -3 (2 a_3 + 8 a_4 x + 20 a_5 x^2 - 3 \gamma y_{x x}) y^{(3)}\\
\text{For $n=8,$ }\quad \mathscr{S} (\mathbf{w}) &= 8\left[3 (a_4 + 5 x (a_5 + 3 a_6 x + 7 a_7 x^2))- 2 \gamma y^{(3)}  \right] y^{(4)}
\end{align*}

This clearly shows that the most general variational symmetry generator for $y^{(n)}=0$ is given by
\begin{equation} \label{vsym40}
\mathbf{w}=  \alpha F_n + \beta G_n +\sum_{k=0} ^{(n-2)/2} a_k V_k
\end{equation}
for $n=4,6, 8.$ \par

Denote by $\mathscr{C} (\mathbf{v})= Q \Delta$ and $\mathscr{F}(\mathbf{v})$  the conservation law and the first integral associated with a given symmetry $\mathbf{v}$ of an {\sc ode} $\Delta_n[y]=0$, where as usual $\mathbf{v}^Q= Q \pd_y.$ Since $V_k^Q= x^k \pd_y$ for all $k$  it follows that
\begin{subequations} \label{w1C&F}
\begin{align}
\mathscr{C}(\mathbf{w}_1)&= \big(P_n(x)\big) \cdot \Delta_n [y],\quad  P_n(x)\equiv P_n(x, (a_k))= \sum_{k=0}^{n-1} a_k x^k   \label{w1C&F1} \\
\mathscr{F}(\mathbf{w}_1) &= \int \mathscr{C}(\mathbf{w}_1) dx +C  \label{w1C&F2}
\end{align}
\end{subequations}
where  $\mathbf{w}_1 =  \sum_{k=0}^{n-1} a_k V_k$ as already noted, and  $C$ is an arbitrary constant of integration, but which we shall omit in the sequel. The integrals in \eqref{w1C&F2} can also be fully evaluated. First of all we note that for each solution symmetry $V_k$ we have
\begin{align}
\mathscr{F} (V_k) &= \sum_{j=0}^k  (-1)^j \left( D^j  x^k\right) y^{(n-1-j)} \notag\\
                  &= \sum_{j=0}^k  (-1)^j  j! \binom{k}{j} x^{k-j} y^{(n-1-j)},\quad (0 \leq k\leq n-1).\label{firstvk}
\end{align}
An important property of the $n$ first integrals $\mathscr{F} (V_k)$  associated with the solution symmetries is that they are linearly  independent. Thus any other first integral of the equation can be expressed in terms of them. More  generally, by computing $\mathscr{F} (\mathbf{w}_1^n)$ for each $n=2,\dots,8$ and rearranging shows that
\begin{equation} \label{firstw1}
\mathscr{F} (\mathbf{w}_1^n) = \sum_{k=0}^{n-1} (-1)^k \left[   D^k P_n (x)  \right] y^{(n-1-k)}
\end{equation}
%%%
The first integrals associated with the symmetry generators of  the  subalgebra $\mathfrak{s}$ of $\mathfrak{g}_n$ have a more complicated pattern which we have not endeavor  to identify.
It clearly follows from \eqref{vqu=1} that
$$
\mathscr{C}(\mathbf{w}_2^n) = \left[-\alpha y_x + \beta( (n-1) y - 2 x y_x)  + \gamma( (1-n) x y + x^2 y_x) \right]\cdot\Delta_n[y].
$$
We thus have the following explicit expressions for the generic generators $\mathbf{w}_2^n,$ for $n=4,6,8.$
\begin{subequations} \label{firtw2}
\begin{align}
& \notag \\
\begin{split}
\mathscr{F} (\mathbf{w}_2^4) &= -2 \gamma  y_x^2+\left(3 y \gamma +(-\beta +\gamma x ) y_x\right) y_{xx}+\frac{1}{2} (\alpha +x (2 \beta -\gamma x )) y_2^2 \\
&\quad+\left(3 y (\beta
-\gamma x )-(\alpha +x (2 \beta -\gamma x )) y_x\right) y^{(3)}\\
\end{split}
& \\
\begin{split}
\mathscr{F} (\mathbf{w}_2^6) &= \frac{9}{2} \gamma  y_2^2+\left(-8 \gamma  y_x+(\beta -\gamma x ) y_{xx}\right) y^{(3)}-\frac{1}{2} (\alpha +x (2 \beta -\gamma x )) y_3^2\\
&\quad + \left(5
y \gamma -3 (\beta -\gamma x ) y_x+\left(\alpha +2 x \beta -x^2 \gamma \right) y_{xx}\right) y^{(4)}\\
&\quad +\left(5 y (\beta -\gamma x )-(\alpha +x
(2 \beta -\gamma x )) y_x\right) y^{(5)}\\
\end{split}
&\\
\begin{split}
\mathscr{F} (\mathbf{w}_2^8) &= -8 \gamma  y_3^2+\left[15 \gamma  y_{xx}+(\gamma x-\beta  ) y^{(3)}\right] y^{(4)}+\frac{1}{2} (\alpha +x (2 \beta -\gamma x )) y_4^2\\
 &\quad -\left[12\gamma  y_x-3 (\beta -\gamma x ) y_{xx}+\left(\alpha +2 x \beta -x^2 \gamma \right) y^{(3)}\right] y^{(5)}\\
&\quad +\left[7 y \gamma -5 (\beta -\gamma x ) y_x+\left(\alpha +2 x \beta -x^2 \gamma \right) y_{xx}\right] y^{(6)}\\
&\quad +\left[7 y (\beta -\gamma x )-(\alpha +x (2 \beta -\gamma x
)) y_x\right] y^{(7)}
\end{split}
\end{align}
\end{subequations}
Using \eqref{firtw2}, one can find in particular $\mathscr{F} (F_n), \mathscr{F} (G_n), \text{ and } \mathscr{F} (H_n)$ by assigning particular values to the linear coefficients $\alpha, \beta$ and $\gamma$ of  $\mathbf{w}_2.$

\subsection{Odd order equations}
As odd order {\small \sc lode}s of maximal symmetry turn out to be all skew-adjoint,  there are no variational symmetries in the sense of \eqref{vsym} associated with them.
%%%
%%%%%%%%%%%%%%%%%%%%%%%
%%%

Moreover, it follows as usual from Theorem \ref{th:wy} and Theorem \ref{th:solsym2}  that in the odd order case to find the divergence symmetry algebra we only need to consider the  subalgebra $\mathfrak{s}$ of $\mathfrak{g}_n.$  This amounts to choosing the most general test vector $\bw$ for divergence symmetries to be of the form $\bw= \bw_2^n.$ The computation of $\mathscr{D} (\mathbf{v})$ for any symmetry operator $\mathbf{v}$ proceeds exactly in the same way as for equations of even orders. Now, in the expression of $\mathscr{D} (\mathbf{w}_2^n)$ the coefficient of the term of highest order, $n+1,$ is $2(\alpha + 2 \beta x - \gamma x^2)$ and the identical vanishing of this term shows indeed that $\mathfrak{s}$ does not generate any divergence symmetry.\par

By Theorem \ref{th:solsym2}, each solution symmetry $V_k$ generates a divergence symmetry, and the associated first integrals $\mathscr{F} (V_k)$ and $\mathscr{F} (\mathbf{w}_1^n)$  naturally have exactly the same expressions as given in \eqref{firstvk} and \eqref{firstw1}, respectively, but for the corresponding odd values of $n=3,5,7.$ Moreover, contrary to the even order case, the homogeneity symmetry $W_y$ generates a divergence symmetry by Theorem \ref{th:wy}, and the associated first integral is given by
\begin{equation} \label{firstwy}
\mathscr{F} (W_y) = \frac{(-1)^{(n-1)/2}}{2} \left(y^{\left(\frac{n-1}{2}\right)}\right)^2 + \sum_{j=0}^{(n-3)/2} (-1)^j y^{(n-1-j)} y^{(j)},
\end{equation}
 for $n=3,5, 7.$ \par

 The results obtained in this section for equations of low orders can naturally be extended to equations of general orders. We denote by $\mathcal{S}_{div}$ and $\mathcal{S}_{var}$ the Lie algebra of divergence symmetries and variational symmetries, respectively, for a given differential equation.
\begin{conj} \label{cj:yn}
For the linear equation $y^{(n)}=0$ with $n \geq 3$ the following holds.
\begin{enumerate}
\item[(a)] For $n$ odd, $\mathcal{S}_{div}= \langle V_0, \dots,  V_{n-1}, W_y\rangle.$ That is, $\mathcal{S}_{div}= \mathcal{A}_n \dotplus \F W_y.$
\item[(b)] For $n$ even,\\ $\mathcal{S}_{div} = \mathcal{A}_n \dotplus \mathfrak{s},$ and $\mathcal{S}_{var} =  \langle V_0, \dots,  V_{\frac{n-2}{2}}\rangle \dotplus  \langle F_n, G_n\rangle,$ with corresponding Lagrangian function given by \eqref{Ln0}.
\item[(c)] For all orders $n,$ the solution symmetries $V_k$ generate $n$ linearly independent first integrals $\mathscr{F} (V_k)$ given by
$$
\mathscr{F} (V_k) = \sum_{j=0}^k  (-1)^j  j! \binom{k}{j} x^{k-j} y^{(n-1-j)},\quad (0 \leq k\leq n-1)
$$
as in \eqref{firstvk}. Thus the most general form of first integrals generated by the solution symmetries is given by
$$
\mathscr{F} (\mathbf{w}_1^n) = \sum_{k=0}^{n-1} (-1)^k \left[   D^k P_n (x)  \right] y^{(n-1-k)},
$$
as in \eqref{firstw1}.
\item[(d)]For $n$ odd the homogeneity symmetry $W_y$ generates  the quadratic first integral with expression
$$
\mathscr{F} (W_y) = \frac{(-1)^{\frac{n-1}{2}}}{2} \left(y^{\left(\frac{n-1}{2}\right)}\right)^2 + \sum_{j=0}^{(n-3)/2} (-1)^j y^{(n-1-j)} y^{(j)},
$$
as in \eqref{firstwy}.
\end{enumerate}
\end{conj}

\begin{proof}
The results in this conjecture have been established for each order $n$ such that $3 \leq n \leq 8.$ For arbitrary orders, parts of the results follow from Theorem \ref{th:wy} and Theorem \ref{th:solsym2} which have been established regardless of the order of the equation, and the rest of the conjecture can be verified order by order.
\end{proof}

Thus, in particular,  for $n$ odd the most general expression for the first integral of $y^{(n)}=0$ is the sum of $\mathscr{F} (\mathbf{w}_1^n)$ and a scalar multiple of $\mathscr{F} (W_y).$

%%%%%%%%%%%%%%%%%%%%%%%%%%%%%%%%%%%%%%%%%%%%%%%
%%%%%%%%%%%%%%%%%%%%%%%%%%%%%%%%%%%%%%%%%%%%%%%%%%%%%%%%%%%%%%%%%%%%%%%%%%%%%%%%%%

\section{Concluding Remarks} \label{s:conclu}
Expressions found in the paper for equations in canonical form concerning their variational symmetries and first integrals as well as similar related results have been almost all expressed  in closed form and for arbitrary orders. Also, while the validity of Theorem \ref{th:wy} has been established only for linear equations of maximal symmetry of the most general form \eqref{nor1}, Theorem \ref{th:solsym2} has however been shown to be valid for all linear equations of any order, without any restriction to the dimension of the symmetry algebra.\par

The study carried out in this paper for linear equations in canonical form is extended in Part 2 \cite{part2} of this series of two papers  to equations of maximal symmetry of the most general form \eqref{nor1}, as well as to the most general case on nonlinear equations of maximal symmetry.

\end{document}